%% file: History_Ravello._Summer_school.tex
\documentclass[10pt, a4paper, twocolumn]{article} 
\usepackage[english]{babel}

\input{structure.tex} 
\usepackage{url}
\usepackage{fancyhdr}
\title{Fifty Years of Excellence: The Summer School of Mathematical Physics of the National Group of Mathematical Physics} 

\author{\authorstyle{Tommaso Ruggeri\textsuperscript{1}, } 
	\authorstyle{Giuseppe Saccomandi\textsuperscript{2} } 
	\newline\newline 
	\textsuperscript{1}\institution{Department of Mathematics and Alma Mater Research Center
on Applied Mathematics $AM^2$, 40123-I BOLOGNA, Italy}\\ 
	\textsuperscript{2}\institution{Department of Engineering, University of Perugia, 06125 PERUGIA, Italy}\\
}

\date{} 

\begin{document}

\maketitle 

\thispagestyle{firstpage} 

\section{Introduction}
“\textit{Education is not the learning of facts, but the training of the mind to think},” Einstein once remarked.  
From this idea comes the conviction that, in order to truly advance the mathematical sciences, isolated initiatives are not sufficient: what is needed are long-term strategies supported by sustained investment in the training of young researchers.  
For this reason, in the second half of the last century, advanced thematic schools were established worldwide, shaping entire generations of scholars.  

In Germany, at the prestigious Mathematisches Forschungsinstitut Oberwolfach, the Oberwolfach Seminars [\cite{1}] have been held since 1981. Aimed at doctoral students and young post-docs (within ten years of their PhD), their purpose is to introduce participants to frontier topics in mathematics.  

In France, the Centre International de Rencontres Mathématiques (CIRM) [\cite{2}], also founded in 1981 by the French mathematical community, has among its missions:  

\begin{quote}
\textit{to ensure transfer of knowledge towards young researchers and doctoral students}.
\end{quote}

In fact, in France the importance of advanced training had already been recognized long before: it suffices to recall the Les Houches Summer School, founded in 1951, primarily devoted to theoretical physics but with strong mathematical components. This very model inspired NATO in 1957 to launch the program \textit{Trained Manpower for Freedom}~\cite{3}, which gave rise to the celebrated NATO Advanced Study Institutes (ASI), now incorporated into the Science for Peace and Security Programme (SPS).

These are just a few examples. The list is naturally much longer: one might mention the El Escorial Meetings in Spain, the activities of the International Centre for Mathematical Meetings at the University of Cantabria, or the Centre International de Mathématiques Pures et Appliquées (CIMPA) [\cite{4}]. Founded in France in 1978 and today a UNESCO center, CIMPA operates on a global scale, aiming to promote training and mathematical research in developing countries. It is supported not only by France but also by Germany, Spain, Norway, and Switzerland.  

The Italian mathematical community has also made significant investments in the training of young scholars. Beginning in the early 1950s, the Fondazione Centro Internazionale Matematico Estivo (C.I.M.E.), operating under the auspices of the Unione Matematica Italiana, promoted with foresight the establishment of summer courses in higher mathematics [\cite{5}]. Later, in 1971, the Scuola Matematica Interuniversitaria (SMI) was created, distinguished from the start by a specific vocation: to organize advanced courses for doctoral students and young fellows, with the goal of providing essential tools for embarking on research. Both the C.I.M.E. and SMI courses continue to this day, taking place in several locations, including the historic sites of Cortona and Perugia.  

Although the mission of the SMI is [\cite{6}]:  

\begin{quote}
\textit{to provide young researchers with a basic training in Mathematics and its applications in various sectors, including Physics, Computer Science, Economics and Finance, Engineering},
\end{quote}

\noindent the courses specifically devoted to mathematical physics (in the Italian sense of the term) have often been limited to an introduction to partial differential equations, overlooking the richness and variety of topics that make up this discipline, situated at the intersection of mathematics, theoretical physics, and engineering sciences.  

Precisely because of its intrinsic interdisciplinarity, mathematical physics plays a crucial role in the development of applied mathematical sciences in industrial, technological, and social contexts.  

To fill this gap and provide young scholars with a broader and more up-to-date training, in 1976 the Scientific Council of the Gruppo Nazionale di Fisica Matematica (GNFM) — then part of the Consiglio Nazionale delle Ricerche (CNR) — decided to establish a dedicated Summer School, capable of offering a more diversified and international program than that already provided by the SMI.  

The goal was ambitious: to create an intensive and stimulating environment in which young researchers could engage with the leading international developments in Mathematical Physics.  

Thus, founded in 1976, the GNFM Summer School of Mathematical Physics has for fifty years represented a benchmark for advanced training in the field. The first edition took place in Bari, followed by two editions in Catania. Beginning in 1979, the School found its permanent home in Ravello, thanks to the initiative of Salvatore Rionero, who directed it until 2017. Since then, the School has been held without interruption, not even during the difficult years of the pandemic.  

Over the years, the School has successfully combined scientific rigor with interdisciplinary openness, establishing itself as one of the most significant initiatives of the GNFM. Today, the GNFM is an integral part of the Istituto Nazionale di Alta Matematica (INDAM), alongside the other national groups in mathematics, and continues to support this School with conviction, which stands as one of the most prestigious and eagerly awaited events of the Italian scientific community in the field.  

\begin{center}
\includegraphics[width=0.4\linewidth]{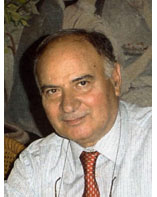}
\captionof{figure}{Professor Salvatore Rionero, who had the insight to bring the School of Mathematical Physics to Ravello in 1979.}
\end{center}

\section{Mission and Objectives}

Since its origins, the School has pursued a dual goal: on the one hand, to provide young researchers with top-level scientific training; on the other, to foster interaction and growth within the scientific community of the field. Over the decades, the School has proved to be a valuable opportunity for at least three main reasons:

\begin{itemize}
\item it has enabled students to access advanced and up-to-date content, often difficult to find, especially in more peripheral universities;
\item it has offered the chance to meet leading figures of the international mathematical community directly;
\item it has created opportunities for networking and scientific exchange among young researchers from different backgrounds and nationalities.
\end{itemize}

Unlike other single-topic schools, the Ravello School offers a wide and articulated range of themes, thus giving participants the chance to engage with subjects beyond their own research area. This openness has often been a source of inspiration and guidance for many young scholars.  

The School takes place every year in September and runs for two weeks, a duration that allows for proper in-depth study of the selected topics. Each day, from Monday to Saturday, includes four hours of lectures in the morning and an afternoon seminar activity, during which the students themselves present and discuss content under the guidance of the lecturers.  

Participants are selected on the basis of their curriculum vitae, with particular attention given to the youngest and most promising profiles.  

Over time, the choice of Ravello has proved strategic for two main reasons. First, it is a relatively secluded location which—while offering an enchanting natural setting—does not present major distractions, thereby fostering concentration and scientific interaction among students. Second, Ravello provides an intimate environment that naturally enhances the sense of community and cohesion among participants. Finally, the beauty of the place serves as a strong attraction for lecturers, who almost always accept with enthusiasm the invitation to teach once contacted by the GNFM Scientific Council. More recently, from a logistical perspective, the use of the prestigious and modern conference hall of Villa Rufolo, provided by the Ravello Foundation, has ensured an excellent setting that enables the School to operate efficiently.  

Several hundreds of young researchers, both Italian and international, have attended the School over the years, and many of them now hold prominent academic positions. The seminar activity assigned to students is not only intended to train them in the synthesis and communication of advanced scientific topics in English, but also often represents a concrete opportunity to be noticed and appreciated by the lecturers. In numerous cases, long-lasting scientific collaborations have emerged from these presentations.  

In this way, the School has significantly contributed to building a solid and international scientific network, connecting the GNFM community with universities and research centers worldwide.  

It should also be noted that the presence of a high-level academic summer school such as the one in Mathematical Physics helps strengthen Ravello’s role as a cultural center. The School can include public lectures or other events involving local citizens, thus creating a bridge between science and the general public. Ravello, therefore, also increases its visibility within the global scientific community, consolidating its reputation as a venue for high-level cultural and academic events.  

\begin{center}
\includegraphics[width=0.4\linewidth]{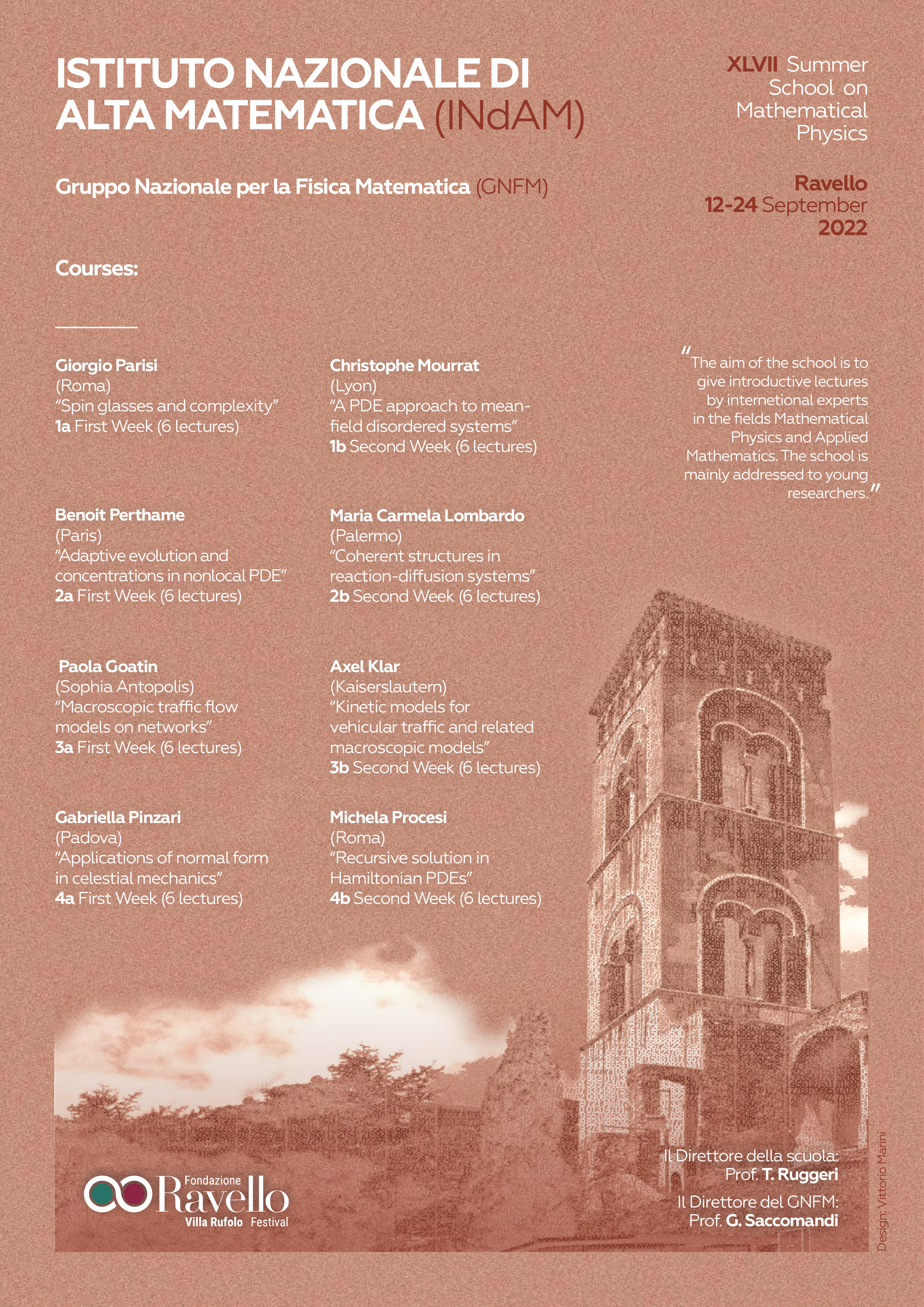}
\includegraphics[width=0.4\linewidth]{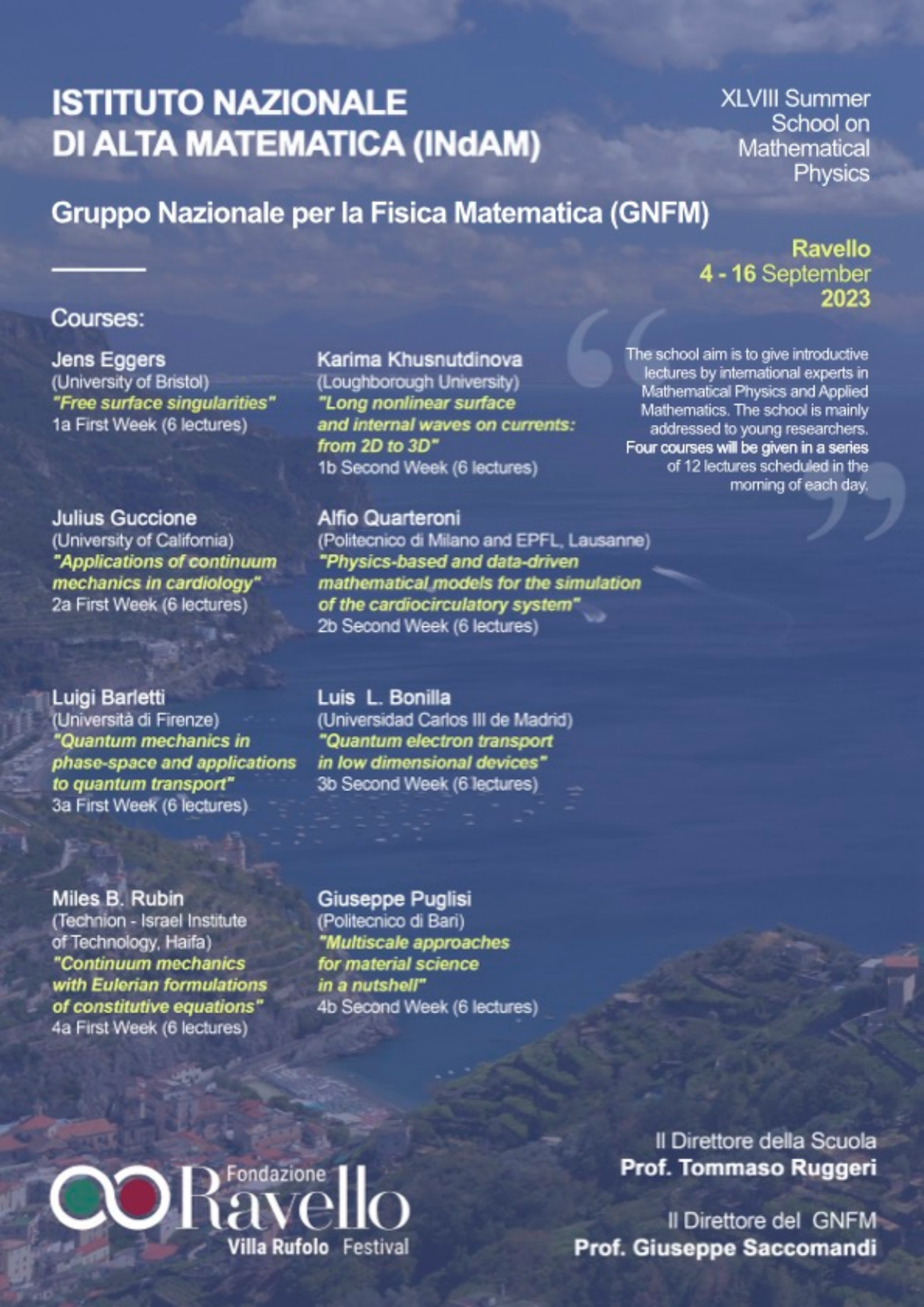}
\captionof{figure}{Posters of the School.}
\end{center}

\section{Topics}

As already mentioned, the Italian notion of Mathematical Physics denotes a discipline far richer and more articulated than its English translation might suggest. In the Italian context, it encompasses a wide range of subjects at the intersection of pure mathematics, theoretical physics, and the applied sciences.  

In Italy, Mathematical Physics has its roots in the tradition of Rational Mechanics, a discipline historically placed within mathematics. Since the 19th century, it has played a central role in scientific education and research, contributing significantly to the development of theoretical models for the description of complex phenomena. The first major research topics concerned the dynamics of rigid bodies, continuum mechanics, potential theory, wave phenomena and, later, relativity and diffusion.  
The Italian approach has always been distinguished by its rigorous mathematical formalism, which made Rational Mechanics a bridge discipline between pure and applied mathematics. Illustrious figures such as Enrico Betti, Tullio Levi-Civita, Antonio Signorini, and Vito Volterra provided fundamental contributions, profoundly influencing both the Italian and international schools.  

From the second half of the 20th century onwards, a profound cultural and scientific transformation took place: professors and researchers in Rational Mechanics increasingly began to identify themselves as Mathematical Physicists, characterized by the systematic use of advanced mathematical tools for the study and formulation of physical models. This shift marked an important turning point in the self-understanding of the discipline, even though in Italy the meaning of “Mathematical Physics” has remained partly distinct from its Anglo-Saxon counterpart, where the emphasis lies primarily on problems of physical interpretation within quantum theory, statistical mechanics, or related fields.  

In recent decades, Italian Mathematical Physics has progressively opened up to new research fields: from statistical mechanics to quantum field theory, from dynamical systems to nonlinear wave theory, and from fluid dynamics to material theory and interdisciplinary applications. At the same time, the evolution of mathematical methods—from functional analysis to partial differential equations and differential geometry—has further expanded the scope and depth of research.  
Today, the Italian Mathematical Physics community maintains a strong identity and a leading international presence, while preserving the historical distinctiveness of its roots in Rational Mechanics. It addresses the challenges of modern and contemporary physics as well as the applications of mathematics to the most diverse fields of society, technology, and the sciences, including biomedical research.  
For these reasons, over the years, the topics addressed in the School have been multiple and diverse, following—and often anticipating—the development of contemporary mathematical physics and contributing to the consolidation of a fruitful dialogue between theoretical formulations, advanced mathematical models, and applications to natural phenomena.  

The main thematic areas addressed in the various editions of the School may be organized as follows:

\begin{itemize}
\item \textbf{Continuum mechanics (fluids and solids).}  
Particular attention has been given to the study of mathematical models describing the behavior of continuous materials, whether fluids or solids. Among the issues addressed are: the dynamics of incompressible and compressible fluids; the mechanics of elastic, viscoelastic, and hyperelastic solids; and the formulation of thermodynamically consistent constitutive models. Considerable attention has also been devoted to the stability of solutions, non-standard boundary conditions, asymptotic regimes, and wave propagation, as well as to the study of complex materials such as liquid crystals.

\item \textbf{Thermodynamics of complex materials.}  
Interest has extended to materials with internal structure, memory effects, or nonlinear behaviors, such as ferromagnetic, electroactive materials, or those subject to coupled thermo-mechanical phenomena. Thermodynamic models have been developed based on variational formulations and entropic principles, with the aim of reconciling microscopic, mesoscopic, and macroscopic descriptions of material behavior.

\item \textbf{Quantum and Statistical Mechanics.}  
The theoretical foundations and mathematical aspects of quantum theories have been explored in depth, both in the canonical and in the functional formulation. Topics covered include the structure of Hilbert spaces, the Schrödinger and Heisenberg equations, the operator algebra formalism, and the connections with statistical mechanics. The latter has been investigated particularly in its deeper aspects: equilibrium states, phase transitions, and emergent phenomena in many-body systems. Moreover, special attention has been devoted to the Boltzmann kinetic theory and its fundamental equation, which describes the evolution of the distribution function and constitutes the bridge between microscopic dynamics and macroscopic transport processes.

\item \textbf{Field theories and geometric structures.}  
The evolution of theoretical physics has led to growing interest in the geometric structures underlying classical and quantum field theories. Tools such as fiber bundles, connections, curvature, and differentiable manifolds have been extensively used to formalize Maxwell’s equations, Yang–Mills theories, and general relativity, as well as to develop approaches to geometric and topological quantization.

\item \textbf{Propagation and transport problems.}  
Another key area has been the study of wave propagation and the transport of mass, energy, or momentum in both continuous and discrete systems. Hyperbolic, parabolic, and dispersive models have been developed and analyzed, with particular focus on nonlinear waves, shock formation, and inverse problems aimed at identifying parameters or sources from observable data.

\item \textbf{Dynamical systems and celestial mechanics.}  
Work in this area has covered various aspects of these disciplines, which occupy a central place in the Italian mathematical physics tradition. Topics have ranged from qualitative to computational analysis, with special emphasis on integrability and symmetry issues, as well as the development of perturbative methods. The research has included both classical issues, such as stereodynamics and the many-body problem, and more modern applications related to vibration mechanics and space mathematics.

\item \textbf{Mathematical analysis applied to complex physical systems.}  
Mathematical analysis has been systematically employed to study physical systems characterized by nonlinear dynamics, emergent behavior, and complex interactions among components. Tools such as differential equations, dynamical systems theory, chaos theory, perturbative methods, and numerical approaches have been used to model phenomena in physical, biological, and technological contexts, with particular attention to bifurcations, multiscale problems, and stability.
\end{itemize}

As already mentioned, this classification is not exhaustive. Over the years, the School has also dealt with emerging and cross-disciplinary topics, such as biomechanical modeling and mathematics for artificial intelligence, reflecting a constant evolution and openness toward the most current scientific challenges.  

Taken as a whole, these topics demonstrate the profound interdisciplinarity and the high conceptual content of mathematical physics, which continues to prove itself a crucial discipline for the understanding and modeling of a wide range of natural, technological, and industrial phenomena.
\begin{center}
\includegraphics[width=0.7\linewidth]{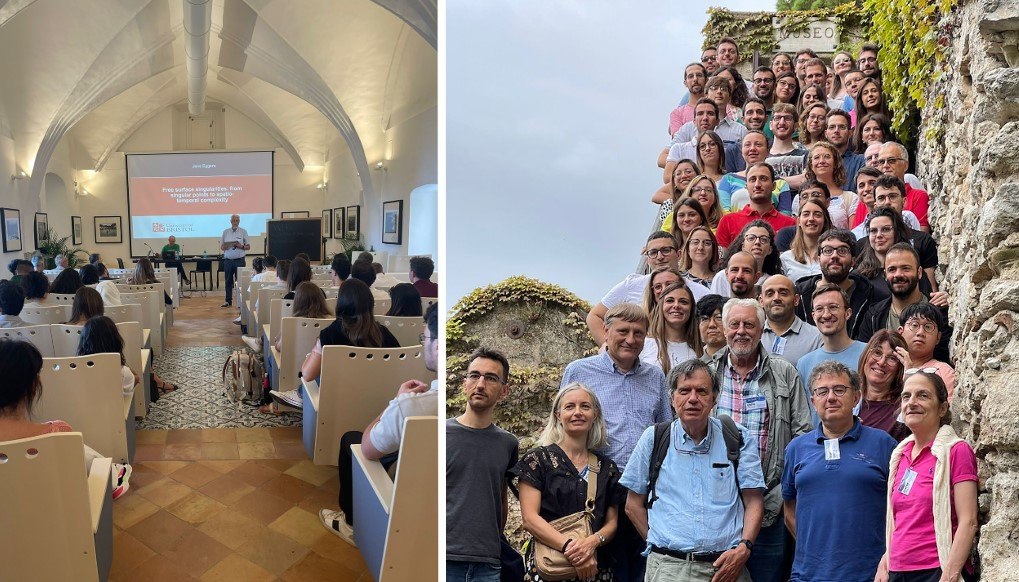}
\captionof{figure}{The lecture hall at Villa Rufolo and a group photo of students and lecturers with Professor Parisi}
\end{center}

\section{Distinguished Lecturers and Representative Courses}

The faculty of the School are, naturally, all eminent scientists.  
Scientific excellence represents, in fact, an essential prerequisite for an initiative of this level.

It should be emphasized, however, that throughout its history the School has had the honor of hosting lecturers of worldwide renown, true celebrities of contemporary science. Just to mention a few names, by way of example:

\begin{itemize}
\item Luis Caffarelli, Abel Prize 2023, recognized for his pioneering work on the regularity of solutions to elliptic and parabolic equations. In 2007, Caffarelli delivered the course: \emph{Nonlinear problems involving fractional diffusions}.
\item Carlo Cercignani (1939–2010), renowned for his fundamental contributions to kinetic theory of gases and the mathematical study of the Boltzmann equation. In 1979, Cercignani delivered the course: \emph{Teoria cinetica dei Gas}.
\item Yvonne Choquet-Bruhat (1923–2025), the first woman elected to the French Académie des Sciences, for her fundamental contributions to the mathematical formulation of general relativity. She taught two courses at Ravello: in 1986, \emph{Supergravity}, and in 2001, \emph{Cauchy Problems in General Relativity}.
\item Peter Lax (1926–2025), Abel Prize 2005, for his fundamental contributions to hyperbolic differential equations and computational mathematics. In 2002, Lax delivered the course: \emph{Non-Linear Wave Propagation}.
\item Pierre-Louis Lions, Fields Medal 1994, for his contributions to the study of differential equations. In 1993, Lions taught the course: \emph{Hartree-Fock and related quantum models of N-body systems}.
\item Jerrold Marsden (1942–2010), who gave essential contributions to differential geometry applied to fluid mechanics and dynamical systems. Winner of the 1973 Gravitational Research Foundation Prize. In 2004, Marsden delivered the course: \emph{Geometrical aspects in Mathematical Physics}.
\item Giorgio Parisi, Nobel Prize in Physics 2021, for his revolutionary discoveries in disordered and complex systems. In 2022, Parisi taught the course: \emph{Spin glasses and complexity}.
\end{itemize}

To provide the reader with a more complete and concrete picture of the School’s activities, Appendix~A presents a selection of particularly significant and representative courses, chosen for their scientific relevance, the topicality of the subjects addressed, or the authority of the lecturers involved\footnote{The full list of courses and lecturers can be found at
\url{https://www.altamatematica.it/gnfm/attivita/scuola-estiva-di-fisica-matematica/}}.

Although partial, this list highlights not only the variety of topics covered, the outstanding scientific level of the lecturers, and their international provenance, but also the School’s ability to bring together scholars from different scientific communities: a clear demonstration of the openness and cultural breadth of the Italian mathematical physics community.

Over the years, in fact, the School has benefited not only from eminent mathematical physicists, but also from analysts such as Gaetano Fichera (1922–1996), numerical analysts such as Alfio Quarteroni, theoretical physicists such as Carlo Rovelli — also known as an essayist and popular science writer — and engineers such as Piero Villaggio (1932–2014).

The School has thus proved capable of speaking the language of mathematics applied to the problems of the real world in a universal way and without any disciplinary prejudice.

\section{Conclusions}

For a school of this kind, reaching fifty years of activity is an extraordinary achievement. It is a tangible proof of scientific dedication and organizational capability. It should be remembered that the organization of an advanced training school is not only a scientific endeavor: its realization also requires a significant logistical and managerial commitment.

The Ravello School has managed to maintain its continuity even during the most difficult times, such as the COVID-19 pandemic, when everything else had come to a halt. Yet, even under those circumstances, the School was held safely and without infections, demonstrating the effectiveness of its organizational framework.

There is another fundamental reason that makes the School indispensable for our community: a prestigious and long-standing school becomes a point of reference, helping to build a shared identity, strengthen institutional visibility, and increase both the national and international prestige of Italian mathematical physics.

It is no coincidence that the Ravello School is now beginning to attract funding from different sponsors, such as local institutions, and is arousing growing interest from the international mathematical community. The Pacific Institute for the Mathematical Sciences (PIMS), for example, has expressed its desire to establish a memorandum of understanding. However, the current structure of the School, based on the voluntary work of a small group of faculty members, makes it difficult to handle the inevitable bureaucracy connected with such agreements.

The Ravello School also represents a significant financial commitment for the \textit{Gruppo Nazionale di Fisica Matematica} (GNFM), whose budget has progressively contracted in recent years. Nevertheless, discontinuing this School would be a serious mistake. Its future is closely tied to that of the GNFM. We must have the courage to continue along the path that our history has handed down to us.

The first challenge is to identify new sources of funding. The School not only guarantees high-quality training for young researchers, but also enhances the value of Ravello as a venue for major cultural and scientific events of international relevance. This opens the door to new opportunities for collaboration, visibility, and growth for the entire city. Local institutions have already begun to recognize this, providing an initial, important financial contribution—support that still has ample potential for growth. Moreover, the Ministry of Universities and Research could consider the School among the activities eligible for support through the \textit{Fondo per la valorizzazione dei Dipartimenti e delle Infrastrutture di Ricerca}, since Ravello represents an international hub for young talents and a platform for interdisciplinary scientific dialogue.

Finally, organizational development must also be envisioned. For example, it would be desirable to adopt an automatic lecture recording system, thus creating a scientific archive of great value.

The challenges, therefore, are many. But there are not many alternatives. As Victor Hugo once reminded us:

\begin{quote}
\textit{Young people have always had the same role in history: to build the future.}
\end{quote}

It is hoped that the GNFM will continue to invest in this School with foresight, passion, and a sense of responsibility. The Ravello School represents, in fact, a unique opportunity to strengthen the unity of the Italian mathematical physics community through a coherent and targeted allocation of resources, while at the same time avoiding the dispersion of available funds.

The \textit{Istituto Nazionale di Alta Matematica} should also continue to feel fully involved—as has already been the case—in this far-reaching project. A School like the one in Ravello is not only an initiative of the GNFM or INdAM, but constitutes a heritage of the entire Italian mathematical community.

\vskip0.5cm
\noindent \textbf{Acknowledgments:} We would like to thank the institutions of Ravello and the Province of Salerno—particularly the Municipality, the \textit{Fondazione Cassa di Risparmio Salernitana}, and the University of Salerno—for supporting the School, as well as the Sunland Agency for managing the logistics.  
We also extend our sincere gratitude to the Scientific Committees of the \textit{Gruppo Nazionale di Fisica Matematica} (GNFM), which have succeeded one another over the years, always providing steady support to the School. Although it is not possible to mention all the committee members individually, we wish to recall Vinicio Boffi and Salvatore Rionero, the first directors of the GNFM and the School, respectively. Finally, a special acknowledgment goes to the Presidents of the \textit{Istituto Nazionale di Alta Matematica} (INdAM), currently led by Professor Cristina Trombetti.

\vskip0.5cm
\begin{center}
\includegraphics[width=0.4\linewidth]{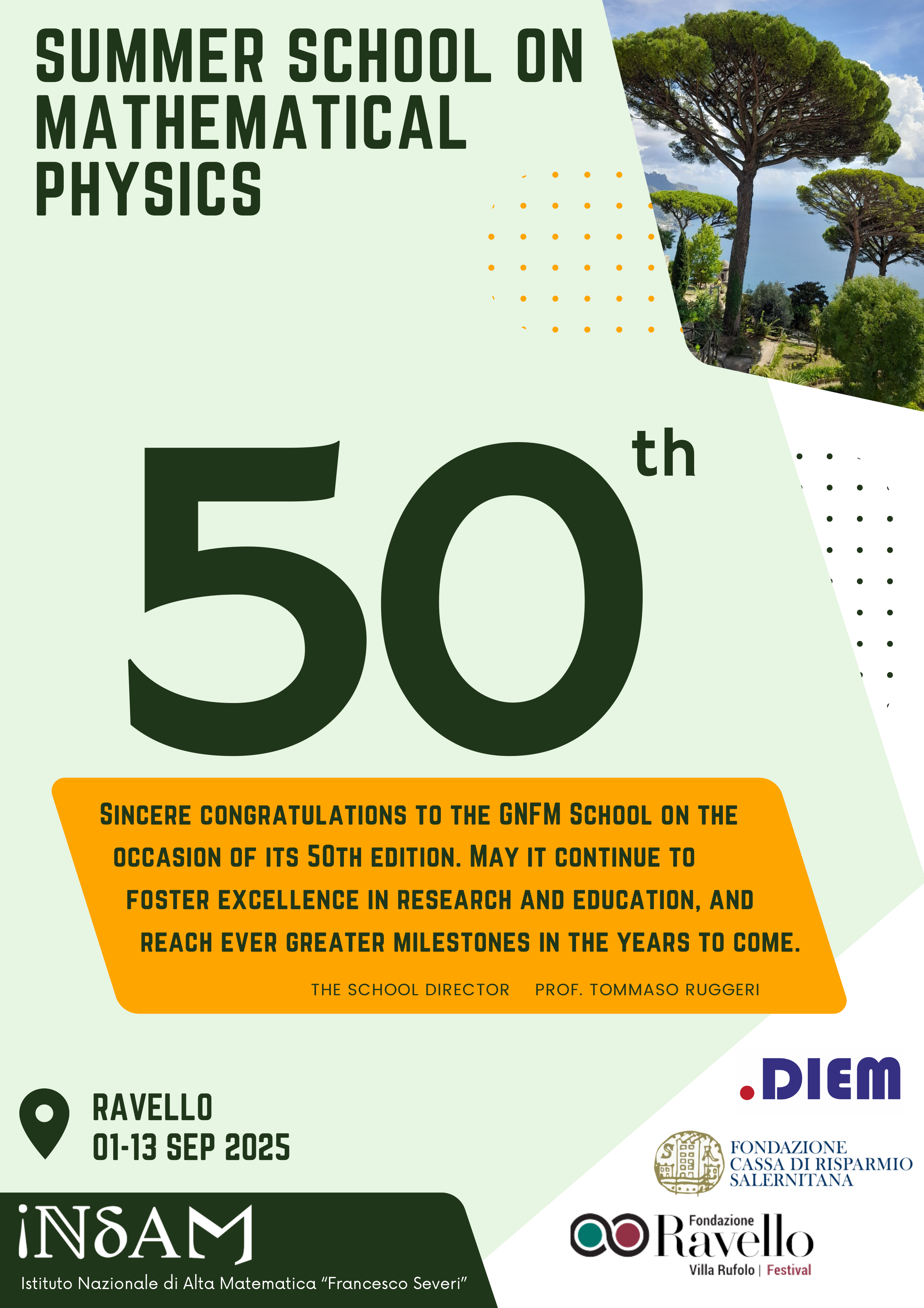}
\captionof{figure}{The 50th-anniversary poster of the School.}
\end{center}

\newpage
\begin{center}
\includegraphics[width=0.7\linewidth]{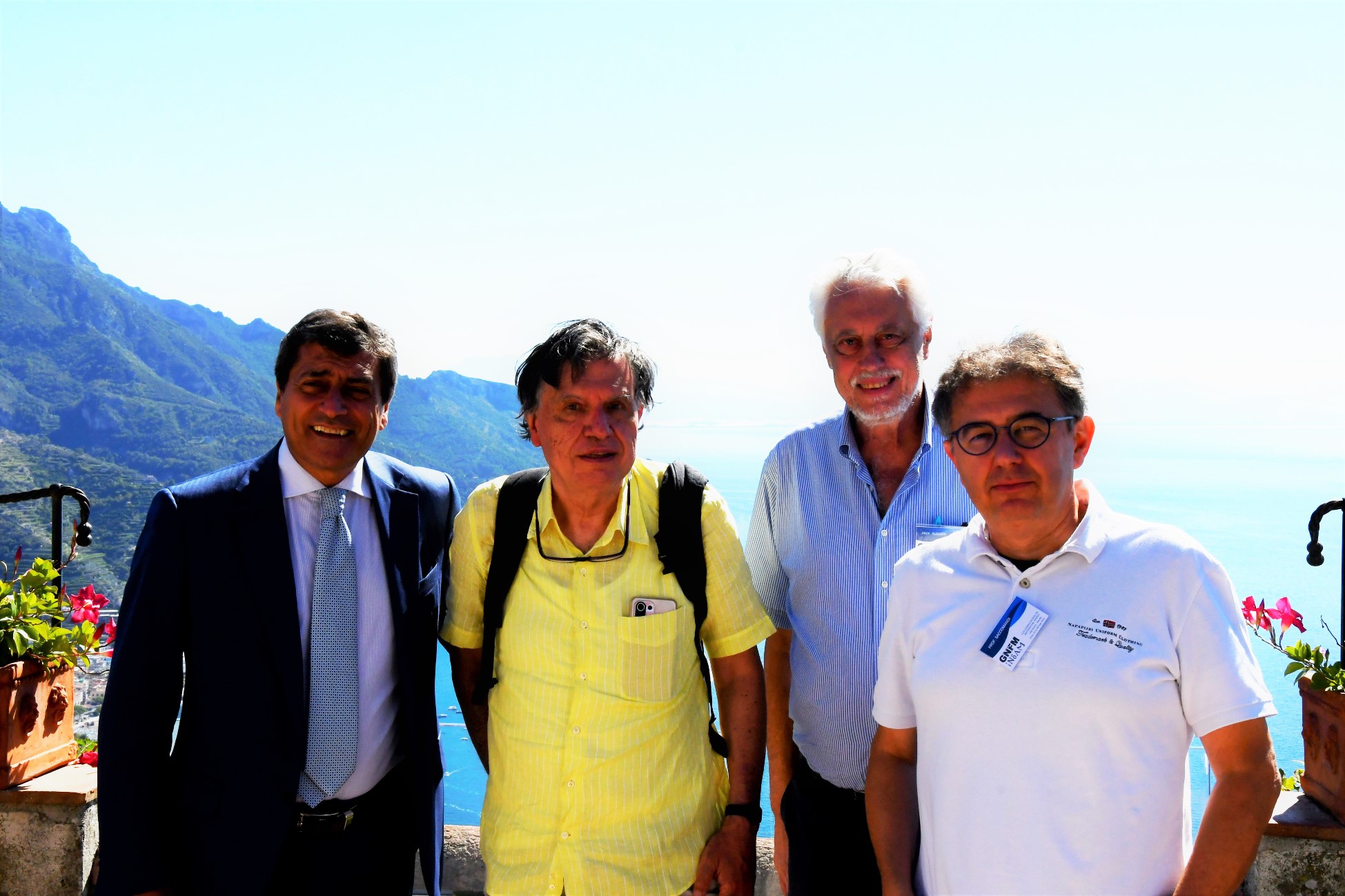}
\captionof{figure}{Professor Parisi with, to his right, the general director of the Villa Rufolo Foundation, Maestro Maurizio Pietrantonio, and to his left Tommaso Ruggeri, current director of the School, and Giuseppe Saccomandi, former director of the GNFM and current member of the INdAM Board.}
\end{center}

\begin{figure}[htbp]
    \centering
    \includegraphics[width=0.7\linewidth]{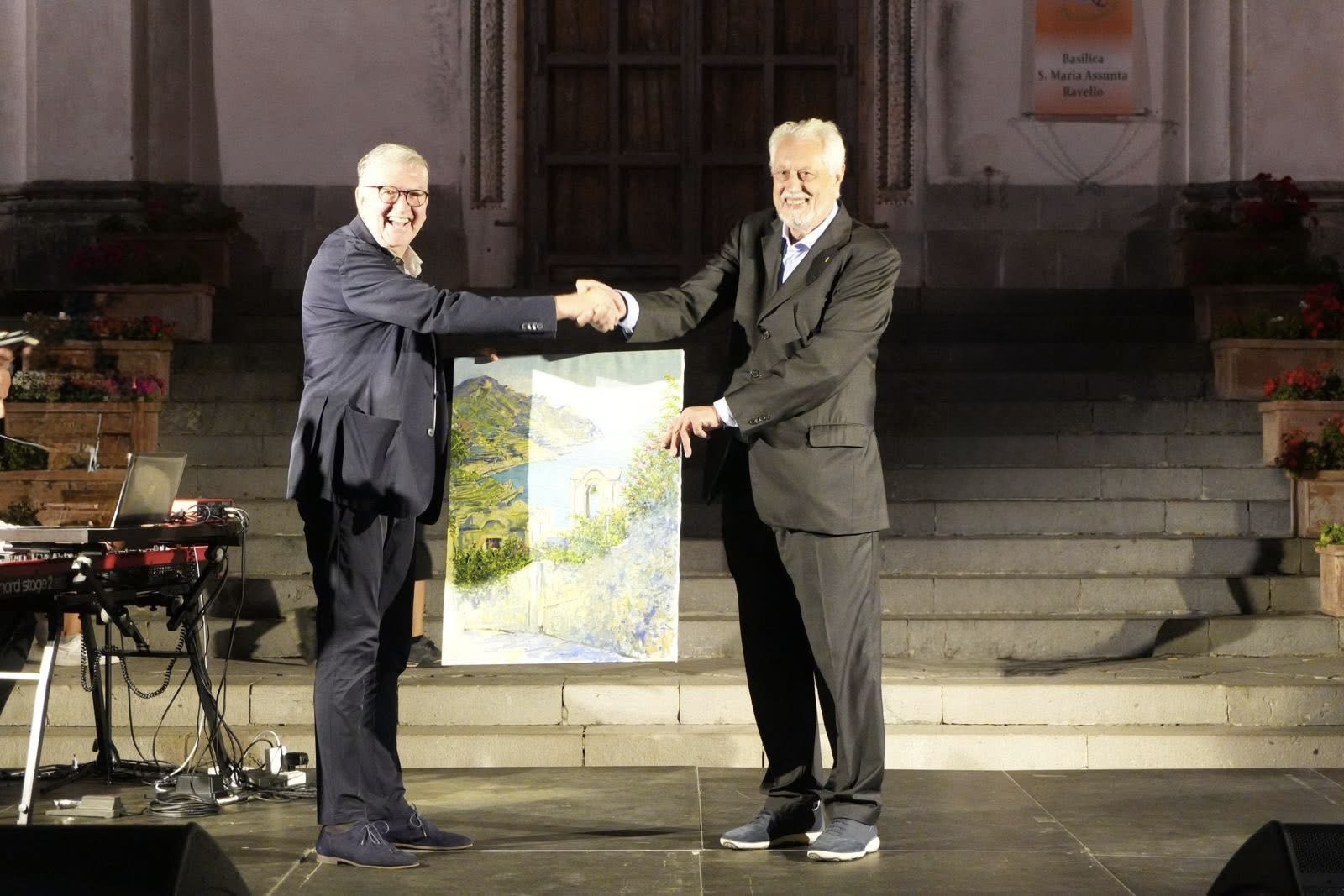}
    \caption{The mayor of Ravello, Paolo Vuilleumier, announced the inclusion of the Summer School of Mathematical Physics, on the occasion of its 50th edition, in the historical register of conferences, presenting the School director with a work by artist Vittorio Abbate.}
    \includegraphics[width=0.7\linewidth]{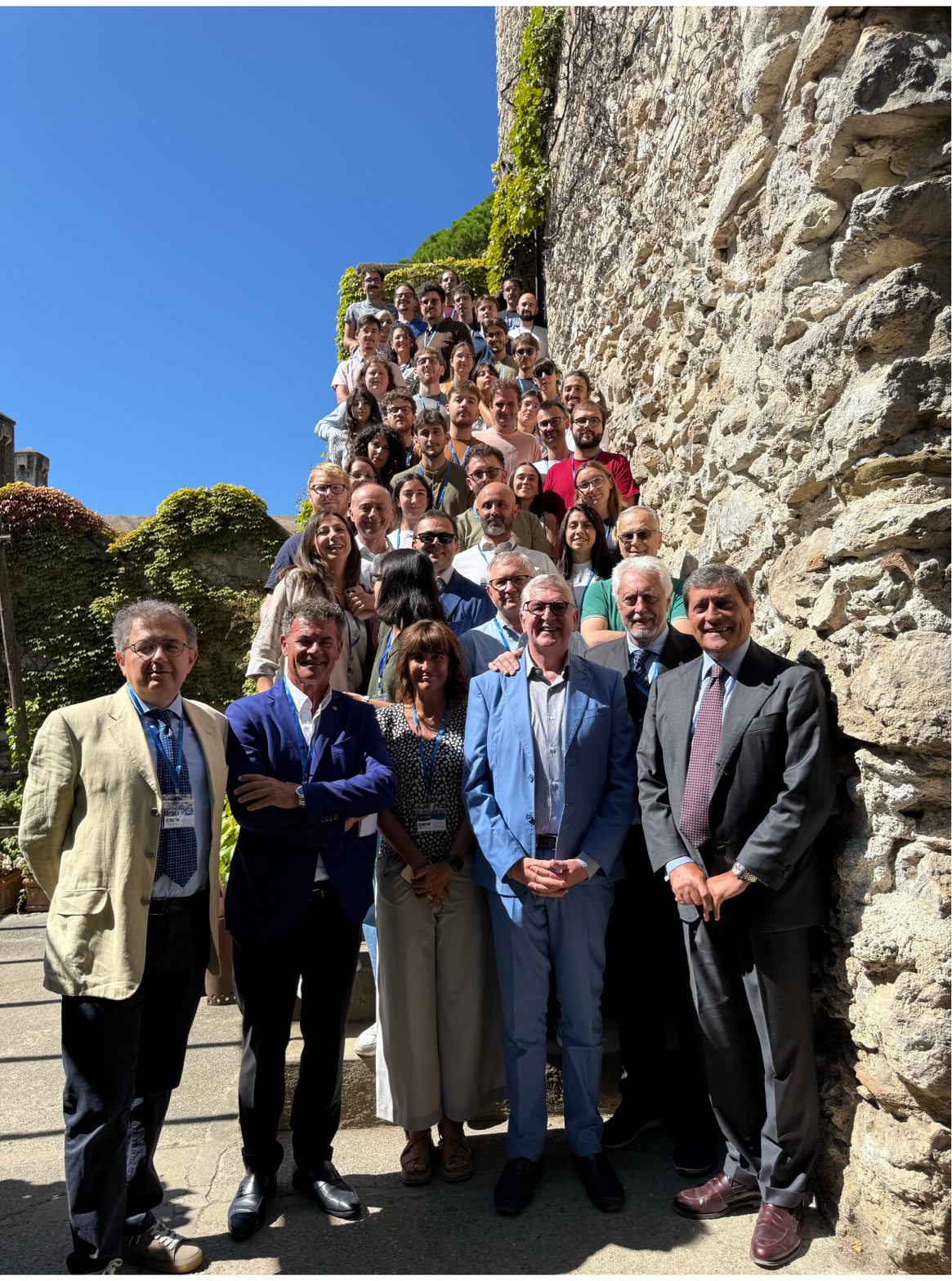}
    \caption{Group photo of the 2025 edition. Front row from left to right: Giuseppe Saccomandi, president of the \textit{Cassa di Risparmio Salernitana} Domenico Credendino, Professor Elisabetta Rocca, mayor of Ravello Paolo Vuilleumier, Tommaso Ruggeri, and director of the Ravello Foundation Maurizio Pietrantonio.}
    \label{fig:ravello}
\end{figure}
\newpage
\noindent \textbf{APPENDIX: Some Notable Courses of the School}
\begin{itemize}
\item[1976] – Meccanica dei continui fluidi – S. Rionero; Thermomechanics – I. Müller; Equazioni alle derivate parziali in problemi di diffusione e di propagazione – M. Primicerio
\item[1977] – Metodi variazionali della fisica matematica – P. Villaggio.
\item[1978] – Relatività ristretta – C. Cattaneo.
\item[1979] – Meccanica delle Vibrazioni – G. Colombo; Teoria cinetica dei gas – C. Cercignani. 
\item[1982] – Stereodinamica – L. Stoppelli.
\item[1984] – Problemi analitici nuovi in fisica matematica classica – G. Fichera.
\item[1985] – Continui relativistici – A. Bressan; Meccanica Stocastica – E. Nelson.
\item[1989] – Phase Transitions – M. Gurtin.
\item[1990] – Mathematical Theory of Multi-phase Fluids – R. Caflish.
\item[1990] – Fluid Dynamics – D. D. Joseph.
\item[1992] – Modelli matematici in biologia – L. Segel. 
\item[1993] – Qualitative methods in nonlinear elasticity – C.M. Dafermos.
\item[1995] –  Lyapounov Methods: Stability for $2-$nd Order Ordinary Differential Systems in $\mathbb{R}^n$, and for Evolutionary Equations in Banach Spaces -- J. Serrin; Topics in Finite Elasticity -- M. Hayes; Topics on PDE’s in Kinetic and Fluid Theory -- B. Perthame. 
\item[1996] -- Aspetti matematici della gravità quantistica -- C. Rovelli.
\item[1997]-- Aspetti matematici della meccanica celeste -- A. Chenciner.
\item[1999] -- Stochastic Differential Equations and Financial Mathematics -- G. Papanicolaou.  
\item[2000] – Measure preserving Maps and Incompressible Fluids – Y. Brenier.
\item[2001] – Fluid Mechanics – V. A. Solonnikov.
\item[2003] – Selected Topics of Kinetic Theory – A. Bobylev.
\item[2005] – Sistemi iperbolici hamiltoniani e loro perturbazioni – B. Dubrovin.
\item[2006] – Nonlinear Electrodynamics and Wave Propagation – G. Boillat.
\item[2007] – Blow up in Continuum Mechanics and scale limits in Statistical Mechanics – E. Presutti.
\item[2009] – Kinetic theory of gases and its applications – K. Aoki.
\item[2013] – Solving Boltzmann Equation, the Green's Function approach – T.-P. Liu.
\item[2015] – Vacuum Dynamics for Compressible Fluids – Z. Xin.
\item[2016] – Introduction to the Mathematical theory of diffusion via PDEs – J.L. Vázquez.
\item[2018] – Analysis of liquid crystals and their defect – J.M. Ball;  Gravitational Waves and Binary Systems – T. Damour.
\item[2024] – Building models of machine learning with mathematical physics – G.M. Rotskoff.
\item[2025] – Interplay between hyperbolic and dispersive equations – Sylvie Benzoni.
\end{itemize}
\newpage 

\end{document}

%% file: structure.tex
\usepackage[english]{babel} 

\usepackage{microtype} 

\usepackage{amsmath,amsfonts,amsthm} 

\usepackage[svgnames]{xcolor} 

\usepackage[hang, small, labelfont=bf, up, textfont=it]{caption} 

\usepackage{booktabs} 

\usepackage{lastpage} 

\usepackage{graphicx} 

\usepackage{enumitem} 
\setlist{noitemsep} 

\usepackage{sectsty} 
\allsectionsfont{\usefont{OT1}{phv}{b}{n}} 

\usepackage{geometry} 

\usepackage[T1]{fontenc} 
\usepackage[utf8]{inputenc} 

\usepackage{XCharter} 

\usepackage{fancyhdr} 
\fancypagestyle{firstpage}{ 
	\fancyhf{}
}

\newcommand{\authorstyle}[1]{{\large\usefont{OT1}{phv}{b}{n}\color{DarkRed}#1}} 

\newcommand{\institution}[1]{{\footnotesize\usefont{OT1}{phv}{m}{sl}\color{Black}#1}} 

\usepackage{titling} 

\newcommand{\HorRule}{\color{DarkGoldenrod}\rule{\linewidth}{1pt}} 

\pretitle{
	\vspace{-30pt} 
	\HorRule\vspace{10pt} 
	\fontsize{32}{36}\usefont{OT1}{phv}{b}{n}\selectfont 
	\color{DarkRed} 
}

\posttitle{\par\vskip 15pt} 

\preauthor{} 

\postauthor{ 
	\vspace{10pt} 
	\par\HorRule 
	\vspace{20pt} 
}

\usepackage{lettrine} 
\usepackage{fix-cm}	

\usepackage{xstring} 
